\documentclass[12pt]{article}
\usepackage{amsfonts}
\usepackage{amssymb}
\usepackage{myart}

\normalbaselines
\input tcilatex
\begin{document}

\title{Different representations of Euclidean geometry and their application
to the space-time geometry}
\author{Yuri A.Rylov}
\date{Institute for Problems in Mechanics, Russian Academy of Sciences,\\
101-1, Vernadskii Ave., Moscow, 119526, Russia.\\
e-mail: rylov@ipmnet.ru\\
Web site: {$http://rsfq1.physics.sunysb.edu/\symbol{126}rylov/yrylov.htm$}\\
or mirror Web site: {$http://gasdyn-ipm.ipmnet.ru/\symbol{126}%
rylov/yrylov.htm$}}
\maketitle

\begin{abstract}
Three different representation of the proper Euclidean geometry are
considered. They differ in the number of basic elements, from which the
geometrical objects are constructed. In E-representation there are three
basic elements (point, segment, angle) and no additional structures.
V-representation contains two basic elements (point, vector) and additional
structure: linear vector space. In $\sigma $-representation there is only
one basic element and additional structure: world function $\sigma =\rho
^{2}/2$, where $\rho $ is the distance. The concept of distance appears in
all representations. However, as a structure, determining the geometry, the
distance appears only in the $\sigma $-representation. The $\sigma $%
-representation is most appropriate for modification of the proper Euclidean
geometry. Practically any modification of the proper Euclidean geometry
turns it into multivariant geometry, where there are many vectors $\mathbf{Q}%
_{0}\mathbf{Q}_{1},\mathbf{Q}_{0}\mathbf{Q}_{1}^{\prime },...$, which are
equal to the vector $\mathbf{P}_{0}\mathbf{P}_{1}$, but they are not equal
between themselves, in general. Concept of multivariance is very important
in application to the space-time geometry. The real space-time geometry is
multivariant. Multivariance of the space-time geometry is responsible for
quantum effects.
\end{abstract}

\section{Introduction}

The proper Euclidean geometry studies mutual dispositions of geometrical
objects (figures) in the space (in the set $\Omega $ of points). Any
geometrical object $\mathcal{O}$ is a subset $\mathcal{O\subset }\Omega $ of
points. Relations between different objects are relations of equivalence,
when two different objects $\mathcal{O}_{1}$ and $\mathcal{O}_{2}$ are
considered to be equivalent $\left( \mathcal{O}_{1}\mathrm{eqv}\mathcal{O}%
_{2}\right) $. The geometrical object $\mathcal{O}_{1}$ is considered to be
equivalent to the geometrical object $\mathcal{O}_{2}$, if after
corresponding displacement the geometrical object $\mathcal{O}_{1}$
coincides with the geometrical object $\mathcal{O}_{2}$.

The geometrical object is considered to be constructed of basic elements
(blocks). There are at least three representations of Euclidean geometry,
which differ in the number and in the choice of basic elements (primary
concepts).

The first representation (Euclidean representation, or E-representation) of
the Euclidean geometry was obtained by Euclid many years ago. The basic
elements in the E-representation are point, segment and angle. The segment
is a segment of the straight line. It consists of infinite number of points.
The segment is determined uniquely by its end points. The angle is a figure,
formed by two segments provided the end of one segment coincides with the
end of other one. Properties of basic elements are described by a system of
axioms. Any geometrical object may be considered to be some composition of
blocks (\textit{point, segment, angle}). The number of the geometric object
constituents may be infinite. The segments determine distances. The angles
determine the mutual orientation of segments. Comparison of different
geometrical objects (figures) $\mathcal{O}_{1}$ and $\mathcal{O}_{2}$ is
produced by their displacement and superposition. If two figures coincide at
superposition, they are considered to be equivalent (equal). The process of
displacement in itself is not formalized in E-representation. However, the
law of the geometric objects displacement is used at the formulation of the
Euclidean geometry propositions and at their proofs.

The second representation (vector representation, or V-representation) of
the Euclidean geometry contains two basic elements (\textit{point, vector}).
From the viewpoint of E-representation the vector is a directed segment of
straight, determined by two points. One point is the origin of the vector,
another point is the end of the vector. But such a definition of vector
takes place only in E-representation, where the vector is a secondary
(derivative) concept. In V-representation the vector is defined
axiomatically as an element of a linear vector space, where two operations
under vectors are defined. These operations (addition of two vectors and
multiplication of a vector by a real number) formalize the law of
displacement of vectors. Strictly, these operations are simply some formal
operations in the linear vector space. They begin to describe the law of
displacement only after introduction of the scalar product of vectors in the
vector space. After introduction of the scalar product the linear vector
space becomes to be the Euclidean space, and the abstract vector may be
considered as a directed segment of straight, determined by two points.
However, it is only interpretation of a vector in the Euclidean
representation. In the V-representation the vector is a primary object. It
is an element of the linear vector space. Nevertheless interpretation of a
vector as a directed segment of straight is very important at construction
of geometrical objects (figures) from points and vectors.

The E-representation contains three basic elements: \textit{point, segment}
and \textit{angle}. The V-representation contains only two basic elements: 
\textit{point} and \textit{vector}. The angle of the E-representation is
replaced by the linear vector space, which is a structure, describing
interrelation of two vectors (directed segments). The vector has some
properties as an element of the linear vector space. Any geometrical figure
may be constructed of points and vectors. It means that the method of
construction of any figure may be described in terms of points and vectors.
The properties of a vector as an element of the linear vector space admit
one to describe properties of displacement of figures and their compositions.

In V-representation the angle appears to be an derivative element. It is
determined by two vectors (by their lengths and by the scalar product
between these vectors). In the V-representation the angle $\theta $ between
two vectors $\mathbf{P}_{0}\mathbf{P}_{1}$ and $\mathbf{P}_{0}\mathbf{P}_{2}$
is defined by the relation 
\begin{equation}
\left\vert \mathbf{P}_{0}\mathbf{P}_{1}\right\vert \cdot \left\vert \mathbf{P%
}_{0}\mathbf{P}_{2}\right\vert \cos \theta =\left( \mathbf{P}_{0}\mathbf{P}%
_{1}.\mathbf{P}_{0}\mathbf{P}_{2}\right)  \label{a1.1}
\end{equation}%
where $\left( \mathbf{P}_{0}\mathbf{P}_{1}.\mathbf{P}_{0}\mathbf{P}%
_{2}\right) $ is the scalar product of vectors $\mathbf{P}_{0}\mathbf{P}_{1}$
and $\mathbf{P}_{0}\mathbf{P}_{2}$. The quantities $\left\vert \mathbf{P}_{0}%
\mathbf{P}_{1}\right\vert \ $and $\left\vert \mathbf{P}_{0}\mathbf{P}%
_{2}\right\vert $ are their lengths%
\begin{equation}
\left\vert \mathbf{P}_{0}\mathbf{P}_{1}\right\vert ^{2}=\left( \mathbf{P}_{0}%
\mathbf{P}_{1}.\mathbf{P}_{0}\mathbf{P}_{1}\right) ,\qquad \left\vert 
\mathbf{P}_{0}\mathbf{P}_{2}\right\vert ^{2}=\left( \mathbf{P}_{0}\mathbf{P}%
_{2}.\mathbf{P}_{0}\mathbf{P}_{2}\right)  \label{a1.2}
\end{equation}

Thus, transition from the E-representation with three basic elements to the
V-representation with two basic elements is possible, provided the
properties of one basic element (vector) are determined by the fact, that
the vector is an element of the linear vector space with the scalar product,
given on it.

Is it possible a further reduction of the number of basic elements in the
representation of the Euclidean geometry? Yes, it is possible. The
representation (in terms world function, or $\sigma $-representation) of the
Euclidean geometry may contain only one basic element (point), provided
there are some constraints, imposed on any two points of the Euclidean space.

The transition from the E-representation to the V-representation reduces the
number of basic elements. Simultaneously this transition generates a new
structure (linear vector space), which determine properties of new basic
element (vector).

The transition from the V-representation to $\sigma $-representation also
reduces the number of basic elements. Only one basic element (\textit{point}%
) remains. Simultaneously this transition replaces the linear vector space
by a new two-point structure (\textit{world function}), which describes
interrelation of two points instead of two vectors. The world function $%
\sigma $ is defined by the relation 
\begin{equation}
\sigma :\qquad \Omega \times \Omega \rightarrow \mathbb{R},\qquad \sigma
\left( P,Q\right) =\sigma \left( Q,P\right) ,\qquad \sigma \left( P,P\right)
=0,\qquad \forall P,Q\in \Omega  \label{a0.1}
\end{equation}%
The world function $\sigma $ is connected with the distance $\rho $ by the
relation%
\begin{equation}
\sigma \left( P,Q\right) =\frac{1}{2}\rho ^{2}\left( P,Q\right) ,\qquad
\forall P,Q\in \Omega  \label{a0.2}
\end{equation}%
The world function of the proper Euclidean space is constrained by a series
of restrictions (see below).

\label{03b} There is a difference between the structures of the
V-representation and that of $\sigma $-representation. Linear vector space
as a structure of the V-representation reflects the symmetry properties of
the Euclidean geometry. These symmetries of the Euclidean space admit the
motion group of the Euclidean space (translations, rotations). These motion
groups admit one to move blocks without deformations and to construct
geometric objects from blocks. The motion groups admit one also to compare
different geometrical objects, moving them in the space.

The world function $\sigma _{\mathrm{E}}$ as a structure of the proper
Euclidean space reflects all properties of the proper Euclidean space, but
not only its symmetries. The world function $\sigma _{\mathrm{E}}$ describes
the properties of the proper Euclidean geometry completely. Any change of
the world function $\sigma _{\mathrm{E}}$ is a deformation of the proper
Euclidean space, which changes its properties.

It should note in this connection, that there are different points of view
on that, what is a geometry. Well known mathematician Felix Klein supposed
that symmetries of a geometry are the most essential properties of the
geometry, and there exist no geometries without a symmetry. For instance, he
wrote that the Riemannian geometry is not a geometry. It is rather a
geography or topography. Such a viewpoint is characteristic for
mathematicians. Most of them believe, that a geometry is a logical
construction, and there exist no nonaxiomatizable geometries.

Alternative viewpoint is characteristic for physicists, who believe, that a
geometry is a science on a shape of geometrical objects and on their mutual
disposition. At such an approach it is of no importance, whether or not the
geometry has any symmetries and whether or not it is axiomatizable. If a
geometry is science on the disposition of geometrical objects, the geometry
is described completely by the world function $\sigma $. Approach of
physicists seems to be more realistic, because the matter distribution
influences on the space-time geometry, and one cannot be sure, that the
space-time geometry is uniform, and it has some symmetries.\label{03e}

Existence of the $\sigma $-representation, containing only one basic
element, means that all geometrical objects and all relations between them
may be recalculated to the $\sigma $-representation, i.e. expressed in terms
of the world function and only in terms of the world function.

If we want to construct a generalized geometry, we are to modify properties
of the proper Euclidean geometry. It means that we are to modify properties
of basic elements of the proper Euclidean geometry. In the E-representation
we have three basic elements (\textit{point, segment, angle}). Their
properties are connected, because finally the segment and the angle are
simply sets of points. Modification the three basic elements cannot be
independent. It is very difficult to preserve connection between the
modified basic elements of a generalized geometry. Nobody does modify the
proper Euclidean geometry in the E-representation.

In V-representation there are two basic elements of geometry (\textit{point,
vector}), and in some cases the modification of the proper Euclidean
geometry is possible. However, the V-representation contains such a
structure, as the linear vector space. One cannot avoid this structure,
because it is not clear, what structure may be used instead. Modification of
the proper Euclidean geometry in the framework of the linear vector space is
restricted rather strong. (It leads to the pseudo-Euclidean geometries and
to the Riemannian geometries). Besides, it appears, that there exist such
modifications of the proper Euclidean geometry, which are incompatible with
the statement that any geometrical object may be constructed of points and
vectors.

The $\sigma $-representation of the proper Euclidean geometry is most
appropriate for modification, because it contains only one basic element (%
\textit{point}). Any modification of the proper Euclidean geometry is
accompanied by a modification of the structure, associated with the $\sigma $%
-representation. This structure is distance (world function), and any
modification of the proper Euclidean geometry is accompanied by a
modification of world function (distance), and vice versa. The meaning of
distance is quite clear.\ This concept appears in all representations of the
proper Euclidean geometry, but only in the $\sigma $-representation the
distance (world function) plays the role of a structure, determining the
geometry. Modification of the distance means a deformation of the proper
Euclidean geometry.

The V-representation appeared in the nineteenth century, and most of
contemporary mathematicians and physicists use this representation. The $%
\sigma $-representation appeared recently in the end of the twentieth
century. Besides, it appears in implicit form in the papers, devoted to
construction of T-geometry \cite{R90,R01}. In these papers the term "$\sigma 
$-representation" was not mentioned, and the $\sigma $-representation was
considered as an evident possibility of the geometry description in terms of
the world function. Apparently, such a possibility was evident only for the
author of the papers, but not for readers. Now we try to correct our default
and to discuss properties of the $\sigma $-representation of the proper
Euclidean geometry.

\section{$\protect\sigma $-representation of the proper Euclidean \newline
geometry}

Let $\Omega $ be the set of points of the Euclidean space. The distance $%
\rho =\rho \left( P_{0},P_{1}\right) $ between two points $P_{0}$ and $P_{1}$
of the Euclidean space is known as the Euclidean metric. In E-representation
as well as in V-representation we have 
\begin{equation}
\rho ^{2}\left( P_{0},P_{1}\right) =\left\vert \mathbf{P}_{0}\mathbf{P}%
_{1}\right\vert ^{2}=\left( \mathbf{P}_{0}\mathbf{P}_{1}.\mathbf{P}_{0}%
\mathbf{P}_{1}\right) ,\qquad \forall P_{0},P_{1}\in \Omega  \label{a2.1}
\end{equation}%
It is convenient to use the world function $\sigma \left( P_{0},P_{1}\right)
=\frac{1}{2}\rho ^{2}\left( P_{0},P_{1}\right) $ as a main characteristic of
the Euclidean geometry. To approach this, one needs to describe properties
of any vector $\mathbf{P}_{0}\mathbf{P}_{1}$ as an element of the linear
vector space in terms of the world function $\sigma \left(
P_{0},P_{1}\right) $, which is associated with the vector $\mathbf{P}_{0}%
\mathbf{P}_{1}$. In $\sigma $-representation the vector $\mathbf{P}_{0}%
\mathbf{P}_{1}$ is defined as the ordered set of two points $\mathbf{P}_{0}%
\mathbf{P}_{1}=\left\{ P_{0},P_{1}\right\} $.

We can to add two vectors, when the end of one vector is the origin of the
other one%
\begin{equation}
\mathbf{P}_{0}\mathbf{Q}_{1}=\mathbf{P}_{0}\mathbf{Q}_{0}+\mathbf{Q}_{0}%
\mathbf{Q}_{1}\ \ \text{or \ }\mathbf{Q}_{0}\mathbf{Q}_{1}=\mathbf{P}_{0}%
\mathbf{Q}_{1}-\mathbf{P}_{0}\mathbf{Q}_{0},  \label{a2.2}
\end{equation}%
\begin{equation}
\mathbf{P}_{1}\mathbf{Q}_{1}=\mathbf{P}_{0}\mathbf{Q}_{1}-\mathbf{P}_{0}%
\mathbf{P}_{1}  \label{a2.2a}
\end{equation}%
Then according to properties of the scalar product in the Euclidean space we
obtain from the second relation (\ref{a2.2})%
\begin{equation}
\left( \mathbf{P}_{0}\mathbf{P}_{1}.\mathbf{Q}_{0}\mathbf{Q}_{1}\right)
=\left( \mathbf{P}_{0}\mathbf{P}_{1}.\mathbf{P}_{0}\mathbf{Q}_{1}\right)
-\left( \mathbf{P}_{0}\mathbf{P}_{1}.\mathbf{P}_{0}\mathbf{Q}_{0}\right)
\label{a2.3}
\end{equation}%
Besides, according to the properties of the scalar product we have%
\begin{equation}
\left\vert \mathbf{P}_{1}\mathbf{Q}_{1}\right\vert ^{2}=\left( \mathbf{P}_{0}%
\mathbf{Q}_{1}-\mathbf{P}_{0}\mathbf{P}_{1}.\mathbf{P}_{0}\mathbf{Q}_{1}-%
\mathbf{P}_{0}\mathbf{P}_{1}\right) =\left\vert \mathbf{P}_{0}\mathbf{Q}%
_{1}\right\vert ^{2}-2\left( \mathbf{P}_{0}\mathbf{P}_{1}.\mathbf{P}_{0}%
\mathbf{Q}_{1}\right) +\left\vert \mathbf{P}_{0}\mathbf{P}_{1}\right\vert
^{2}  \label{a2.4}
\end{equation}%
Using definition of the world function%
\begin{equation}
\sigma \left( P_{0},P_{1}\right) =\sigma \left( P_{1},P_{0}\right) =\frac{1}{%
2}\left\vert \mathbf{P}_{0}\mathbf{P}_{1}\right\vert ^{2},\qquad \forall
P_{0},P_{1}\in \Omega  \label{a2.5}
\end{equation}%
we obtain from (\ref{a2.4}) and (\ref{a2.5})%
\begin{equation}
\left( \mathbf{P}_{0}\mathbf{P}_{1}.\mathbf{P}_{0}\mathbf{Q}_{1}\right)
=\sigma \left( P_{0},P_{1}\right) +\sigma \left( P_{0},Q_{1}\right) -\sigma
\left( P_{1},Q_{1}\right) ,\qquad \forall P_{0},P_{1},Q_{1}\in \Omega
\label{a2.6}
\end{equation}

It follows from (\ref{a2.3}) and (\ref{a2.6}), that for any two vectors $%
\mathbf{P}_{0}\mathbf{P}_{1}$ and $\mathbf{Q}_{0}\mathbf{Q}_{1}$ the scalar
product has the form%
\begin{equation}
\left( \mathbf{P}_{0}\mathbf{P}_{1}.\mathbf{Q}_{0}\mathbf{Q}_{1}\right)
=\sigma \left( P_{0},Q_{1}\right) +\sigma \left( P_{1},Q_{0}\right) -\sigma
\left( P_{0},Q_{0}\right) -\sigma \left( P_{1},Q_{1}\right) ,\qquad \forall
P_{0},P_{1},Q_{0},Q_{1}\in \Omega  \label{a2.7}
\end{equation}%
Setting $Q_{0}=P_{0}$ in (\ref{a2.7}) and comparing with (\ref{a2.6}), we
obtain%
\begin{equation}
\sigma \left( P_{0},P_{0}\right) =\frac{1}{2}\left\vert \mathbf{P}_{0}%
\mathbf{P}_{0}\right\vert ^{2}=0,\qquad \forall P_{0}\in \Omega
\label{a2.7a}
\end{equation}

$n$ vectors $\mathbf{P}_{0}\mathbf{P}_{1},\mathbf{P}_{0}\mathbf{P}_{2},...\
\ \mathbf{P}_{0}\mathbf{P}_{n}\ $ are linear dependent, if and only if the
Gram's determinant $F_{n}\left( \mathcal{P}^{n}\right) $, $\mathcal{P}%
^{n}\equiv \left\{ P_{0},P_{1},...P_{n}\right\} $ vanishes%
\begin{equation}
F_{n}\left( \mathcal{P}^{n}\right) \equiv \det \left\vert \left\vert \left( 
\mathbf{P}_{0}\mathbf{P}_{i}.\mathbf{P}_{0}\mathbf{P}_{k}\right) \right\vert
\right\vert =0,\qquad i,k=1,2,...n  \label{a2.8}
\end{equation}%
In the $\sigma $-representation the condition (\ref{a2.8}) is written in the
developed form%
\begin{equation}
F_{n}\left( \mathcal{P}^{n}\right) \equiv \det \left\vert \left\vert \sigma
\left( P_{0},P_{i}\right) +\sigma \left( P_{0},P_{k}\right) -\sigma \left(
P_{i},P_{k}\right) \right\vert \right\vert =0,\qquad i,k=1,2,...n
\label{a2.9}
\end{equation}

If $\sigma $ is the world function of $n$-dimensional Euclidean space, it
satisfies the following relations.

I. Definition of the dimension and introduction of the rectilinear
coordinate system: 
\begin{equation}
\exists \mathcal{P}^{n}\equiv \left\{ P_{0},P_{1},...P_{n}\right\} \subset
\Omega ,\qquad F_{n}\left( \mathcal{P}^{n}\right) \neq 0,\qquad F_{k}\left( {%
\Omega }^{k+1}\right) =0,\qquad k>n  \label{g2.5}
\end{equation}%
where $F_{n}\left( \mathcal{P}^{n}\right) $\ is the Gram's determinant (\ref%
{a2.9}). Vectors $P_{0}P_{i}$, $\;i=1,2,...n$\ are basic vectors of the
rectilinear coordinate system $K_{n}$\ with the origin at the point $P_{0}$.
The covariant metric tensor $g_{ik}\left( \mathcal{P}^{n}\right) $, \ $%
i,k=1,2,...n$\ and the contravariant one $g^{ik}\left( \mathcal{P}%
^{n}\right) $, \ $i,k=1,2,...n$\ in $K_{n}$\ are defined by the relations 
\begin{equation}
\sum\limits_{k=1}^{k=n}g^{ik}\left( \mathcal{P}^{n}\right) g_{lk}\left( 
\mathcal{P}^{n}\right) =\delta _{l}^{i},\qquad g_{il}\left( \mathcal{P}%
^{n}\right) =\left( \mathbf{P}_{0}\mathbf{P}_{i}.\mathbf{P}_{0}\mathbf{P}%
_{l}\right) ,\qquad i,l=1,2,...n  \label{a1.5b}
\end{equation}%
\begin{equation}
F_{n}\left( \mathcal{P}^{n}\right) =\det \left\vert \left\vert g_{ik}\left( 
\mathcal{P}^{n}\right) \right\vert \right\vert \neq 0,\qquad i,k=1,2,...n
\label{g2.6}
\end{equation}

II. Linear structure of the Euclidean space: 
\begin{equation}
\sigma \left( P,Q\right) =\frac{1}{2}\sum\limits_{i,k=1}^{i,k=n}g^{ik}\left( 
\mathcal{P}^{n}\right) \left( x_{i}\left( P\right) -x_{i}\left( Q\right)
\right) \left( x_{k}\left( P\right) -x_{k}\left( Q\right) \right) ,\qquad
\forall P,Q\in \Omega  \label{a1.5a}
\end{equation}%
where coordinates $x_{i}\left( P\right) ,$\ $i=1,2,...n$\ of the point $P$\
are covariant coordinates of the vector $\mathbf{P}_{0}\mathbf{P}$, defined
by the relation\textit{\ } 
\begin{equation}
x_{i}\left( P\right) =\left( \mathbf{P}_{0}\mathbf{P}_{i}.\mathbf{P}_{0}%
\mathbf{P}\right) ,\qquad i=1,2,...n  \label{b12}
\end{equation}

III: The metric tensor matrix $g_{lk}\left( \mathcal{P}^{n}\right) $\ has
only positive eigenvalues 
\begin{equation}
g_{k}>0,\qquad k=1,2,...,n  \label{a15c}
\end{equation}

IV. The continuity condition: the system of equations 
\begin{equation}
\left( \mathbf{P}_{0}\mathbf{P}_{i}.\mathbf{P}_{0}\mathbf{P}\right)
=y_{i}\in \mathbb{R},\qquad i=1,2,...n  \label{b14}
\end{equation}%
considered to be equations for determination of the point $P$\ as a function
of coordinates $y=\left\{ y_{i}\right\} $,\ \ $i=1,2,...n$\ has always one
and only one solution.\textit{\ }Conditions I $\div $ IV contain a reference
to the dimension $n$\ of the Euclidean space.

One can show that conditions I $\div $ IV are the necessary and sufficient
conditions of the fact that the set $\Omega $ together with the world
function $\sigma $, given on $\Omega $, describes the $n$-dimensional
Euclidean space \cite{R90}.

Thus, in the $\sigma $-representation the Euclidean geometry contains only
one primary geometrical object: the point. Any two points are described by
the world function $\sigma $, which satisfies conditions I $\div $ IV. Any
geometrical figure and any relation can be described in terms of the world
function and only in terms of the world function.

In the $\sigma $-representation the vector $\mathbf{P}_{0}\mathbf{P}_{1}$ is
a ordered set $\left\{ P_{0},P_{1}\right\} $ of two points. Scalar product
of two vectors is defined by the relation (\ref{a2.7}).

Two vectors $\mathbf{P}_{0}\mathbf{P}_{1}$ and $\mathbf{Q}_{0}\mathbf{Q}_{1}$
are collinear (linear dependent), if%
\begin{equation}
\mathbf{P}_{0}\mathbf{P}_{1}\parallel \mathbf{Q}_{0}\mathbf{Q}_{1}:\qquad
\left( \mathbf{P}_{0}\mathbf{P}_{1}.\mathbf{Q}_{0}\mathbf{Q}_{1}\right)
^{2}=\left\vert \mathbf{P}_{0}\mathbf{P}_{1}\right\vert ^{2}\left\vert 
\mathbf{Q}_{0}\mathbf{Q}_{1}\right\vert ^{2}  \label{a2.14}
\end{equation}%
Two vectors $\mathbf{P}_{0}\mathbf{P}_{1}$ and $\mathbf{Q}_{0}\mathbf{Q}_{1}$
are equivalent (equal), if%
\begin{equation}
\mathbf{P}_{0}\mathbf{P}_{1}=\mathbf{Q}_{0}\mathbf{Q}_{1}:\qquad \left( 
\mathbf{P}_{0}\mathbf{P}_{1}.\mathbf{Q}_{0}\mathbf{Q}_{1}\right) =\left\vert 
\mathbf{P}_{0}\mathbf{P}_{1}\right\vert \cdot \left\vert \mathbf{Q}_{0}%
\mathbf{Q}_{1}\right\vert \wedge \left\vert \mathbf{P}_{0}\mathbf{P}%
_{1}\right\vert =\left\vert \mathbf{Q}_{0}\mathbf{Q}_{1}\right\vert
\label{a2.15}
\end{equation}

Vector $\mathbf{S}_{0}\mathbf{S}_{1}$ with the origin at the given point $%
S_{0}$ is the sum of two vectors $\mathbf{P}_{0}\mathbf{P}_{1}$ and $\mathbf{%
Q}_{0}\mathbf{Q}_{1}$,%
\begin{equation}
\mathbf{S}_{0}\mathbf{S}_{1}=\mathbf{S}_{0}\mathbf{R}+\mathbf{RS}_{1},
\label{a2.16}
\end{equation}%
if the points $S_{1}$ and $R$ satisfy the relations%
\begin{equation}
\mathbf{S}_{0}\mathbf{R=P}_{0}\mathbf{P}_{1},\qquad \mathbf{RS}_{1}=\mathbf{Q%
}_{0}\mathbf{Q}_{1}  \label{a2.17}
\end{equation}%
In the developed form it means that the points $S_{1}$ and $R$ satisfy the
relations 
\begin{eqnarray}
\left( \mathbf{S}_{0}\mathbf{R.P}_{0}\mathbf{P}_{1}\right) &=&\left\vert 
\mathbf{S}_{0}\mathbf{R}\right\vert \cdot \left\vert \mathbf{P}_{0}\mathbf{P}%
_{1}\right\vert ,\qquad \left\vert \mathbf{S}_{0}\mathbf{R}\right\vert
=\left\vert \mathbf{P}_{0}\mathbf{P}_{1}\right\vert  \label{a2.18} \\
\left( \mathbf{RS}_{1}.\mathbf{Q}_{0}\mathbf{Q}_{1}\right) &=&\left\vert 
\mathbf{RS}_{1}\right\vert \cdot \left\vert \mathbf{Q}_{0}\mathbf{Q}%
_{1}\right\vert ,\qquad \left\vert \mathbf{RS}_{1}\right\vert =\left\vert 
\mathbf{Q}_{0}\mathbf{Q}_{1}\right\vert  \label{a2.19}
\end{eqnarray}%
where scalar products are expressed via corresponding world functions by the
relation (\ref{a2.7}). The points $P_{0},P_{1},Q_{0},Q_{1},S_{0}$ are given.
One can determine the point $R$ from two equations (\ref{a2.18}). As far as
the world function satisfy the conditions I $\div $ IV, the geometry is the
Euclidean one, and the equations (\ref{a2.18}) have one and only one
solution for the point $R$. When the point $R$ has been determined, one can
determine the point $S_{1}$, solving two equations (\ref{a2.19}). They also
have one and only one solution for the point $S_{1}$.

Result of summation does not depend on the choice of the origin $S_{0}$ in
the following sense. Let the sum of vectors $\mathbf{P}_{0}\mathbf{P}_{1}$
and $\mathbf{Q}_{0}\mathbf{Q}_{1}$ be defined with the origin at the point $%
S_{0}^{\prime }$ by means of the conditions%
\begin{equation}
\mathbf{S}_{0}^{\prime }\mathbf{S}_{1}^{\prime }=\mathbf{S}_{0}^{\prime }%
\mathbf{R}^{\prime }+\mathbf{R}^{\prime }\mathbf{S}_{1}^{\prime }
\label{a2.20}
\end{equation}%
where the points $S_{1}^{\prime }$ and $R^{\prime }$ satisfy the relations%
\begin{equation}
\mathbf{S}_{0}^{\prime }\mathbf{R}^{\prime }\mathbf{=P}_{0}\mathbf{P}%
_{1},\qquad \mathbf{R}^{\prime }\mathbf{S}_{1}^{\prime }=\mathbf{Q}_{0}%
\mathbf{Q}_{1}  \label{a2.21}
\end{equation}%
Then in force of conditions I $\div $ IV the geometry is the Euclidean one,
and there is one and only one solution of equations (\ref{a2.21}) and%
\begin{equation}
\mathbf{S}_{0}\mathbf{S}_{1}=\mathbf{S}_{0}^{\prime }\mathbf{S}_{1}^{\prime }
\label{a2.22}
\end{equation}%
in the sense, that%
\begin{equation}
\mathbf{S}_{0}\mathbf{S}_{1}=\mathbf{S}_{0}^{\prime }\mathbf{S}_{1}^{\prime
}:\qquad \left( \mathbf{S}_{0}\mathbf{S}_{1}.\mathbf{S}_{0}^{\prime }\mathbf{%
S}_{1}^{\prime }\right) =\left\vert \mathbf{S}_{0}\mathbf{S}_{1}\right\vert
\cdot \left\vert \mathbf{S}_{0}^{\prime }\mathbf{S}_{1}^{\prime }\right\vert
\wedge \left\vert \mathbf{S}_{0}\mathbf{S}_{1}\right\vert =\left\vert 
\mathbf{S}_{0}^{\prime }\mathbf{S}_{1}^{\prime }\right\vert  \label{a2.23}
\end{equation}

Let us stress that in the $\sigma $-representation the sum of two vectors
does not depend on the choice of the origin of the resulting vector, because
the world function satisfies the Euclideaness conditions I $\div $ IV. If
the world function does not satisfy the conditions I $\div $ IV, the result
of summation may depend on the origin $S_{0}$, as well as on the order of
vectors $\mathbf{P}_{0}\mathbf{P}_{1}$, $\mathbf{Q}_{0}\mathbf{Q}_{1}$ at
summation. Besides, the result may be multivariant even for a fixed point $%
S_{0}$, because the solution of equations (\ref{a2.17}) may be not unique.

Multiplication of the vector $\mathbf{P}_{0}\mathbf{P}_{1}$ by a real number 
$\alpha $ gives the vector $\mathbf{S}_{0}\mathbf{S}_{1}$ with the origin at
the point $S_{0}$ 
\begin{equation}
\mathbf{S}_{0}\mathbf{S}_{1}=\alpha \mathbf{P}_{0}\mathbf{P}_{1}
\label{a2.24}
\end{equation}%
Here the points $P_{0},P_{1},S_{0}$ are given, and the point $S_{1}$ is
determined by the relations%
\begin{equation}
\left( \mathbf{S}_{0}\mathbf{S}_{1}.\mathbf{P}_{0}\mathbf{P}_{1}\right) =%
\text{sgn}\left( \alpha \right) \left\vert \mathbf{S}_{0}\mathbf{S}%
_{1}\right\vert \cdot \left\vert \mathbf{P}_{0}\mathbf{P}_{1}\right\vert
,\qquad \left\vert \mathbf{S}_{0}\mathbf{S}_{1}\right\vert =\left\vert
\alpha \right\vert \cdot \left\vert \mathbf{P}_{0}\mathbf{P}_{1}\right\vert
\label{a2.25}
\end{equation}%
Because of conditions I $\div $ IV there is one and only one solution of
equations (\ref{a2.25}) and the solution does not depend on the point $S_{0}$%
.

\section{Uniqueness of operations in the proper \newline
Euclidean geometry}

To define operations under vectors (equality, summation, multiplication) in
the $\sigma $-representation, one needs to solve algebraic equations (\ref%
{a2.15}), (\ref{a2.21}), (\ref{a2.24}), which are reduced finally to
equations (\ref{a2.15}), defining the equality operation.

Equality of two vectors $\mathbf{P}_{0}\mathbf{P}_{1}$ and $\mathbf{Q}_{0}%
\mathbf{Q}_{1}$, $\left( \mathbf{P}_{0}\mathbf{P}_{1}=\mathbf{Q}_{0}\mathbf{Q%
}_{1}\right) $ is defined by two algebraic equations%
\begin{equation}
\left( \mathbf{P}_{0}\mathbf{P}_{1}.\mathbf{Q}_{0}\mathbf{Q}_{1}\right)
=\left\vert \mathbf{P}_{0}\mathbf{P}_{1}\right\vert \cdot \left\vert \mathbf{%
Q}_{0}\mathbf{Q}_{1}\right\vert ,\qquad \left\vert \mathbf{P}_{0}\mathbf{P}%
_{1}\right\vert =\left\vert \mathbf{Q}_{0}\mathbf{Q}_{1}\right\vert
\label{a3.1}
\end{equation}%
The number of equations does not depend on the dimension of the proper
Euclidean space.

In V-representation the number of equations, determining the equality of
vectors, is equal to the dimension of the Euclidean space. To define
equality of $\mathbf{P}_{0}\mathbf{P}_{1}$ and $\mathbf{Q}_{0}\mathbf{Q}_{1}$%
, one introduces the rectilinear coordinate system $K_{n}$ with basic
vectors $\mathbf{e}_{0},\mathbf{e}_{1},...,\mathbf{e}_{n-1}$ and the origin
at the point $O$. Covariant coordinates $x_{k}=\left( \mathbf{P}_{0}\mathbf{P%
}_{1}\right) _{k}$ and $y_{k}=\left( \mathbf{Q}_{0}\mathbf{Q}_{1}\right)
_{k} $ are defined by the relations 
\begin{equation}
x_{k}\equiv \left( \mathbf{P}_{0}\mathbf{P}_{1}\right) _{k}=\left( \mathbf{e}%
_{k}.\mathbf{P}_{0}\mathbf{P}_{1}\right) ,\qquad y_{k}\equiv \left( \mathbf{Q%
}_{0}\mathbf{Q}_{1}\right) _{k}=\left( \mathbf{e}_{k}.\mathbf{Q}_{0}\mathbf{Q%
}_{1}\right) \qquad k=0,1,...n-1  \label{a3.2}
\end{equation}%
The equality equations of vectors $\mathbf{P}_{0}\mathbf{P}_{1}$ and $%
\mathbf{Q}_{0}\mathbf{Q}_{1}$ in V-representation have the form%
\begin{equation}
x_{k}=y_{k},\qquad k=0,1,...,n-1  \label{a3.3}
\end{equation}

According to the linear structure of the Euclidean space (\ref{a1.5a}) and
due to definition of the scalar product (\ref{a2.7}) we obtain%
\begin{equation}
\left\vert \mathbf{P}_{0}\mathbf{P}_{1}\right\vert
^{2}=g^{ik}x_{i}x_{k},\qquad \left\vert \mathbf{Q}_{0}\mathbf{Q}%
_{1}\right\vert ^2=g^{ik}y_{i}y_{k},\qquad \left( \mathbf{P}_{0}\mathbf{P}%
_{1}.\mathbf{Q}_{0}\mathbf{Q}_{1}\right) =g^{ik}x_{i}y_{k}  \label{a3.4}
\end{equation}%
where%
\begin{equation}
g^{ik}g_{lk}=\delta _{l}^{i},\qquad g_{il}=\left( \mathbf{e}_{i}.\mathbf{e}%
_{l}\right) ,\qquad i,k=0,1,...,n-1  \label{a3.5}
\end{equation}%
A summation $0\div \left( n-1\right) $ is produced on repeated indices.

Due to relations (\ref{a3.4}) two equality relations (\ref{a3.1}) take the
form%
\begin{equation}
g^{kl}x_{k}y_{l}=g^{kl}x_{k}x_{l},\qquad g^{kl}x_{k}x_{l}=g^{kl}y_{k}y_{l}
\label{a3.6}
\end{equation}%
By means of the second equation (\ref{a3.6}) the first equation (\ref{a3.6})
may be written in the form%
\begin{equation}
g^{kl}\left( x_{k}-y_{k}\right) \left( x_{l}-y_{l}\right) =0  \label{a3.7}
\end{equation}%
According to the III condition (\ref{a15c}) the matrix of the metric tensor $%
g^{kl}$ has only positive eigenvalues. In this case the equation (\ref{a3.7}%
) has only trivial solution for $x_{k}-y_{k}$, $k=0,1,...n-1$. Then the
equation (\ref{a3.7}) is equivalent to $n$ equations%
\begin{equation}
x_{l}=y_{l},\qquad l=0,1,...,n-1  \label{a3.8}
\end{equation}

Thus, in the $\sigma $-representation two equations (\ref{a3.1}) of the
vector equality are equivalent to $n$ equations (\ref{a3.3}) of vector
equality in V-representation. This equivalency is conditioned by the linear
structure (\ref{a1.5a}) of the Euclidean space and by the positive
distinctness of the Euclidean metric. In the pseudo-Euclidean space the
matrix of the metric tensor has eigenvalues of different sign. In this case
the relations (\ref{a3.7}) and (\ref{a3.8}) cease to be equivalent.

\section{Generalization of the proper Euclidean \newline
geometry}

Let us consider a simple example of the proper Euclidean geometry
modification. The matrix $g^{ik}$ of the metric tensor in the rectilinear
coordinate system $K_{n}$ with basic vectors $\mathbf{e}_{0},\mathbf{e}%
_{1},...\mathbf{e}_{n-1}$ is modified in such a way that its eigenvalues
have different signs. We obtain pseudo-Euclidean geometry. For simplicity we
set the dimension $n=4$ and eigenvalues $g_{0}=1,$ $g_{1}=g_{2}=g_{3}=-1$.
It is the well known space of Minkowski (pseudo-Euclidean space of index $1$%
). In the space of Minkowski the definition (\ref{a3.1}) of two vectors
equality in $\sigma $-representation does not coincide, in general, with the
two vectors equality (\ref{a3.3}) in V-representation. The two definitions
are equivalent for timelike vectors ($g^{kl}x_{k}x_{l}>0$), and they are not
equivalent for spacelike vectors ($g^{kl}x_{k}x_{l}<0$).

Indeed, using definition (\ref{a3.1}), we obtain the relation (\ref{a3.7})
for any vectors of the same length. It means that the vector with
coordinates $x_{k}-y_{k}$ is an isotropic vector. However the sum, as well
as the difference of two timelike vectors of the same length in the
Minkowski space is either a timelike vector, or a zeroth vector. Hence, the
relations (\ref{a3.8}) take place. Thus, definitions of equality (\ref{a3.6}%
) and (\ref{a3.3}) coincide for timelike vectors.

In the case of spacelike vectors $x_{k}$ and $y_{k}$ their difference $%
x_{k}-y_{k}$ may be an isotropic vector. The relation (\ref{a3.7}) states
this. It means that equality of two spacelike vectors of the same length in
the $\sigma $-representation of the Minkowski space is multivariant, i.e.
there are many spacelike vectors $\mathbf{Q}_{0}\mathbf{Q}_{1},\mathbf{Q}_{0}%
\mathbf{Q}_{1}^{\prime },...$ which are equivalent to the spacelike vector $%
\mathbf{P}_{0}\mathbf{P}_{1}$, but the vectors $\mathbf{Q}_{0}\mathbf{Q}_{1},%
\mathbf{Q}_{0}\mathbf{Q}_{1}^{\prime },...$ are not equivalent between
themselves, in general.

The result is rather unexpected. Firstly, the definition of equality appears
to be different in V-representation and in $\sigma $-representation.
Secondly, the equality definition in $\sigma $-representation appears to be
multivariant, what is very unusual.

In the V-representation the equality of two vectors is single-variant by
definition of a vector as an element of the linear vector space. In the $%
\sigma $-representation the equality of two vectors is defined by the
world-function. The pseudo-Euclidean space is a result of a deformation of
the proper Euclidean space, when the Euclidean world function is replaced by
the pseudo-Euclidean world function. The definition of the vector equality
via the world function remains to be the equality definition in all
Euclidean spaces.

For brevity we shall use different names for geometries with different
definition of the vectors equality. Let geometry with the vectors equality
definition (\ref{a3.3}), (\ref{a3.2}) be the Minkowskian geometry, whereas
the geometry with the vectors equality definition (\ref{a3.1}) will be
referred to as the $\sigma $-Minkowskian geometry. The world function is the
same in both geometries. The geometries differ in the structure of the
linear vector space and, in particular, in definition of the vector
equality. Strictly, there is no linear vector space in the $\sigma $%
-Minkowskian geometry. The linear vector space is not necessary for
formulation of the Euclidean geometry, as well as for formulation of the $%
\sigma $-Minkowskian geometry. The Minkowskian geometry and the $\sigma $%
-Minkowskian one are constructed in different representations. The geometry,
constructed on the basis of the world function is a more general geometry,
because the world function exists in both geometries.

\label{01b}\textit{Important remark. }In application of the space-time
geometry of Minkowski to the space-time the spacelike world lines and
spacelike vectors are not used. The difference between the geometry of
Minkowski and $\sigma $-minkowskian geometry appears to be unessential from
this viewpoint. However, this remark does not concern more general
space-time geometries.\label{01e}

What of the two equality definitions (\ref{a2.15}), or (\ref{a3.3}) are
valid? Maybe, both? Apparently, both geometries are possible as abstract
constructions. The Minkowskian geometry and the $\sigma $-Minkowskian
geometry differ in such a property as the multivariance of spacelike vectors
equality. In the geometry, constructed on the basis of the linear vector
space, the multivariance is absent in principle (by definition). On the
contrary, in the $\sigma $-representation the equivalent vector is
determined as a solution of the equations (\ref{a2.15}). For arbitrary world
function one can guarantee neither existence, nor uniqueness of the
solution. They may be guaranteed, only if the world function satisfies the
conditions I $\div $ IV. It means, that the multivariance is a general
property of the generalized geometry, whereas the single-variance of the
proper Euclidean geometry is a special property of the proper Euclidean
geometry. The special properties of the proper Euclidean geometry are
described by the conditions I $\div $IV. All these conditions contain a
reference to the dimension of the Euclidean space. The single-variance of
the proper Euclidean geometry is a special property, which is determined by
the form of the world function, but it does not contain a reference to the
dimension, because it is valid for proper Euclidean space of any dimension.

From formal viewpoint the $\sigma $-Minkowskian geometry is more consistent,
because the definition (\ref{a3.1}) does not depend on the choice of the
coordinate system. In the Minkowskian geometry the two vectors equality
definition contains a reference to the coordinate system, which may be
considered as an additional structure introduced to the geometry of
Minkowski. After this the geometry if Minkowski should be qualified as a
fortified geometry. It means that a physical geometry is equipped by some
additional geometric structure. This structure suppresses multivariance of
the spacelike vectors equality.

We may dislike this fact, because we habituated ourselves to single-variance
of the two vectors equivalence. However, we may not ignore the fact, that
the multivariance is a natural property of geometry. It means that we are to
use the two vector equivalence in the form (\ref{a2.15}).

Let us discuss corollaries of the new definition of the vector equivalence
for the real space-time. In the Euclidean space the straight line, passing
through the points $P_{0},P_{1}$ is defined by the relation 
\begin{equation}
\mathcal{T}_{P_{0}P_{1}}=\left\{ R|\mathbf{P}_{0}\mathbf{R}\parallel \mathbf{%
P}_{0}\mathbf{P}_{1}\right\}  \label{a4.3}
\end{equation}%
In the $\sigma $-representation the parallelism $\mathbf{P}_{0}\mathbf{R}%
\parallel \mathbf{P}_{0}\mathbf{P}_{1}$ of vectors $\mathbf{P}_{0}\mathbf{R}$
and $\mathbf{P}_{0}\mathbf{P}_{1}$ is defined by the relation (\ref{a2.14})
in the $\sigma $-Minkowskian space-time. In the Minkowskian space-time the
vector parallelism $\mathbf{P}_{0}\mathbf{R}\parallel \mathbf{P}_{0}\mathbf{P%
}_{1}$ is defined by four equations, describing proportionality of
components of vectors $\mathbf{P}_{0}\mathbf{R}$ and $\mathbf{P}_{0}\mathbf{P%
}_{1}$. In the $\sigma $-Minkowskian space-time the straight (\ref{a4.3}) is
a one-dimensional line for the timelike vector $\mathbf{P}_{0}\mathbf{P}_{1}$%
, whereas the straight (\ref{a4.3}) is a three-dimensional surface (two
planes) for the spacelike vector $\mathbf{P}_{0}\mathbf{P}_{1}$. In the
Minkowskian space-time all straights (timelike and spacelike) are
one-dimensional lines. Thus, the straights, generated by the spacelike
vector, are described differently in the Minkowskian geometry and in the $%
\sigma $-Minkowskian one.

The timelike straight in the geometry of Minkowski describes the world line
of a free particle, whereas the spacelike straight is believed to describe a
hypothetical particle taxyon. The taxyon is not yet discovered. In the $%
\sigma $-Minkowskian space-time geometry this fact is explained as follows.
The taxyon, if it exists, is described by two isotropic three-dimensional
planes, which is an envelope to a set of light cones, having its vertices on
some one-dimensional straight line. Such a taxyon has not been discovered,
because one looks for it in the form of the one-dimensional straight line.
In the conventional Minkowskian space-time the taxyon has not been
discovered, because there exist no particles, moving with the speed more
than the speed of the light. Thus, the absence of taxyon is explained on the
geometrical level in the $\sigma $-Minkowskian space-time geometry, whereas
the absence of taxyon is explained on the dynamic level in the conventional
Minkowskian space-time geometry.

Let us consider the non-Riemannian space-time geometry $\mathcal{G}_{\mathrm{%
d}}$, described by the world function $\sigma _{\mathrm{d}}$%
\begin{equation}
\sigma _{\mathrm{d}}=\sigma _{\mathrm{M}}+\mathrm{sgn}\left( \sigma _{%
\mathrm{M}}\right) d,\qquad d=\lambda _{0}^{2}=\frac{\hbar }{2bc}=\text{const%
}  \label{a4.4}
\end{equation}%
\label{04b}%
\begin{equation}
\mathrm{sgn}\left( x\right) =\left\{ 
\begin{array}{l}
1,\ \ \text{if}\ \ x>0 \\ 
0,\ \ \ \ \text{if\ \ }x=0 \\ 
-1,\ \ \text{if}\ \ x<0%
\end{array}%
\right. ,  \label{a4.4a}
\end{equation}%
where $\sigma _{\mathrm{M}}$ is the world function of the geometry of
Minkowski, $\hbar $ is the quantum constant, $c$ is the speed of the light
and $b$ is some universal constant $\left[ b\right] =$g/cm. The constant $%
\lambda _{0}$ is some universal length.

The space-time geometry (\ref{a4.4}) is discrete , because in this
space-time geometry there are no vectors $\mathbf{P}_{0}\mathbf{P}_{1}$,
whose length $\left\vert \mathbf{P}_{0}\mathbf{P}_{1}\right\vert $ be small
enough, i.e. 
\begin{equation}
\left\vert \mathbf{P}_{0}\mathbf{P}_{1}\right\vert ^{4}\notin \left(
0,4\lambda _{0}^{4}\right) ,\qquad \forall P_{0},P_{1}\in \Omega
\label{a4.4b}
\end{equation}%
In other words, the space-time geometry (\ref{a4.4}) has no close points.
The space-time geometry, where there are no close points should be qualified
as a discrete geometry. It is rather unexpected that a discrete space-time
geometry may be uniform and isotropic as the geometry (\ref{a4.4}). It is
unexpected, that a discrete geometry may be given on a continuous manifold.
But this puzzle is connected with the fact, that usually one uses
V-representation, where a discrete geometry is given on a lattice space.
Another unexpected properties of a discrete space-time geometry can be found
in \cite{R2008}. \label{04e}

In the space-time geometry (\ref{a4.4}) a free pointlike particle motion is
described by a chain $\mathcal{T}_{\mathrm{br}}$ of connected segments $%
\mathcal{T}_{\left[ P_{k}P_{k+1}\right] }$ of straight 
\begin{equation}
\mathcal{T}_{\mathrm{br}}=\dbigcup\limits_{k}\mathcal{T}_{\left[ P_{k}P_{k+1}%
\right] },  \label{a4.5}
\end{equation}%
\begin{equation}
\mathcal{T}_{\left[ P_{k}P_{k+1}\right] }=\left\{ R|\sqrt{2\sigma _{\mathrm{d%
}}\left( P_{k},R\right) }+\sqrt{2\sigma _{\mathrm{d}}\left( P_{k+1},R\right) 
}=\sqrt{2\sigma _{\mathrm{d}}\left( P_{k},P_{k+1}\right) }\right\}
\label{a4.6}
\end{equation}%
The particle 4-momentum $\mathbf{p}$ is described by the vector $\mathbf{P}%
_{k}\mathbf{P}_{k+1}$ 
\begin{equation}
\mathbf{p}=bc\mathbf{P}_{k}\mathbf{P}_{k+1},\qquad m=b\left\vert \mathbf{P}%
_{k}\mathbf{P}_{k+1}\right\vert =b\mu ,\qquad k=0,\pm 1,...  \label{a4.7}
\end{equation}%
Here $m$ is the particle mass, $\mu $ is the geometrical particle mass, i.e.
the particle mass expressed in units of length. Description of $\mathcal{G}_{%
\mathrm{d}}$ is produced in the $\sigma $-representation, because
description of non-Riemannian geometry in the V-representation is
impossible. In the geometry $\mathcal{G}_{\mathrm{d}}$ segments $\mathcal{T}%
_{\left[ P_{k}P_{k+1}\right] }$ of the timelike straight are multivariant in
the sense, that $\mathcal{T}_{\left[ P_{k}P_{k+1}\right] }$ is a
cigar-shaped three-dimensional surface, but not a one-dimensional segment.

For the free particle the adjacent links (4-momenta) $\mathcal{T}_{\left[
P_{k}P_{k+1}\right] }$ and $\mathcal{T}_{\left[ P_{k+1}P_{k+2}\right] }$ are
equivalent in the sense of (\ref{a2.15})%
\begin{eqnarray}
\mathbf{P}_{k}\mathbf{P}_{k+1}\text{eqv}\mathbf{P}_{k+1}\mathbf{P}_{k+2} &:&
\nonumber \\
\left( \mathbf{P}_{k}\mathbf{P}_{k+1}.\mathbf{P}_{k+1}\mathbf{P}%
_{k+2}\right) &=&\left\vert \mathbf{P}_{k}\mathbf{P}_{k+1}\right\vert \cdot
\left\vert \mathbf{P}_{k+1}\mathbf{P}_{k+2}\right\vert ,\qquad \left\vert 
\mathbf{P}_{k}\mathbf{P}_{k+1}\right\vert =\left\vert \mathbf{P}_{k+1}%
\mathbf{P}_{k+2}\right\vert  \label{a4.8a}
\end{eqnarray}%
This equivalence is multivariant in the sense that at fixed link $\mathbf{P}%
_{k}\mathbf{P}_{k+1}$ the adjacent link $\mathbf{P}_{k+1}\mathbf{P}_{k+2}$
wobbles with the characteristic angle $\theta =\sqrt{d/\mu ^{2}}=\sqrt{%
b\hbar /m^{2}c}$. The shape of the chain with wobbling links appears to be
random. Statistical description of the random world chain appears to be
equivalent to the quantum description in terms of the Schr\"{o}dinger
equation \cite{R91}. Thus, quantum effects are a corollary of the
multivariance. Here the multivariance is taken into account on the level of
the space-time geometry. In conventional quantum theory the space-time
geometry is single-variant, whereas the dynamics is multivariant, because,
when one replaces the conventional dynamic variable by a matrix or by
operator, one introduces the multivariance in dynamics.

Thus, the multivariance is a natural property of the real world \label{05b}%
(especially of microcosm).\label{05e} At the conventional approach one
avoids the multivariance "by hand" from geometry and introduces it "by hand"
in dynamics, to explain quantum effects. It would be more reasonable to
remain the multivariance in the space-time geometry, because it appears
there naturally. Besides, multivariance of the space-time geometry has
another corollaries, other than quantum effects. In particular,
multivariance of the space-time geometry is a reason of the restricted
divisibility of real bodies into parts (atomism) \cite{R2006}.

\section{Multivariance and possibility of the geometry axiomatization}

In the E-representation of the proper Euclidean geometry one supposes that
any geometrical object can be constructed of basic elements (points,
segments, angles). In the V-representation any geometrical object is
supposed to be constructed of basic elements (points, vectors). Deduction of
propositions of the proper Euclidean geometry from the system of axioms
reproduces the process of the geometrical object construction. There are
different ways of the basic elements application for a construction of some
geometrical object, because such basic elements as segments and vectors are
simply sets of points. In the analogical way a proposition of the proper
Euclidean geometry may be obtained by different proofs, based on the system
of axioms.

\label{06b}Let us consider a simple example of the three-dimensional proper
Euclidean space. In the Cartesian coordinates $\mathbf{x}=\left\{
x,y,z\right\} $ the world function has the form 
\begin{equation}
\sigma _{\mathrm{E}}\left( \mathbf{x},\mathbf{x}^{\prime }\right) =\frac{1}{2%
}\left( \left( x-x^{\prime }\right) ^{2}+\left( y-y^{\prime }\right)
^{2}+\left( z-z^{\prime }\right) ^{2}\right)   \label{a5.1}
\end{equation}%
We consider a ball $\mathcal{B}$ with the boundary%
\begin{equation}
x^{2}+y^{2}+z^{2}=R^{2}  \label{a5.2}
\end{equation}%
where $R$ is the radius of the ball. The ball $\mathcal{B}$ may by
considered as constructed only of the set $\mathcal{S}$ of the straight
segments $\mathcal{T}\left( -x,y,z;x,y,z\right) $, $x^{2}+y^{2}+z^{2}=R^{2}$%
\begin{equation}
\mathcal{B}=\dbigcup\limits_{x^{2}+y^{2}+z^{2}=R^{2}}\mathcal{T}\left(
-x,y,z;x,y,z\right)   \label{a5.4}
\end{equation}%
Any segment $\mathcal{T}\left( -x,y,z;x,y,z\right) $ is a segment of the
length $\sqrt{R^{2}-y^{2}-z^{2}}$. Its center is placed at the point $%
\left\{ 0,y,z\right\} $. The segment $\mathcal{T}\left( x,y,z;-x,y,z\right) $
is the set of points $P^{\prime }$ with coordinates $\left( x^{\prime
},y,z\right) $%
\begin{equation}
\mathcal{T}\left( -x,y,z;x,y,z\right) =\left\{ \left\{ x^{\prime }|x^{\prime
2}\leq x^{2}\right\} ,y,z\right\}   \label{a5.3}
\end{equation}%
Different segments $\mathcal{T}\left( -x,y,z;x,y,z\right) $ have no common
points, and any point $P$ of the ball $\mathcal{B}$ belongs to one and only
one of segments of the set $\mathcal{S}$.\label{06e}

Let us consider the ball $\mathcal{B}$ with the boundary 
\begin{equation}
x^{2}+y^{2}+z^{2}=R^{2}-d/2  \label{a5.3a}
\end{equation}
But now we consider non-Riemannian geometry $\mathcal{G}_{\mathrm{d}}$,
described by the world function $\sigma _{\mathrm{d}}$ 
\begin{equation}
\sigma _{\mathrm{d}}=\left\{ 
\begin{array}{l}
\sigma _{\mathrm{E}}-d,\quad \text{if\quad }\sigma _{\mathrm{E}}\geq 2d \\ 
\frac{1}{2}\sigma _{\mathrm{E}},\quad \text{if\quad }\sigma _{\mathrm{E}}<2d%
\end{array}%
\right. ,\qquad d=\text{const, }d>0  \label{a5.5}
\end{equation}%
where $\sigma _{\mathrm{E}}$ is the proper Euclidean world function (\ref%
{a5.1}). Thus, if $\sigma _{\mathrm{d}}\gg d$ the world function $\sigma _{%
\mathrm{d}}$ distinguishes slightly from $\sigma _{\mathrm{E}}$.

The straight segment $\mathcal{T}_{\left[ P_{0}P_{1}\right] }$ between the
points $P_{0}$ and $P_{1}$ is the set of points 
\begin{equation}
\mathcal{T}_{\left[ P_{0}P_{1}\right] }=\left\{ R|\sqrt{2\sigma _{\mathrm{d}%
}\left( P_{0},R\right) }+\sqrt{2\sigma _{\mathrm{d}}\left( P_{1},R\right) }=%
\sqrt{2\sigma _{\mathrm{d}}\left( P_{0},P_{1}\right) }\right\}  \label{a5.6}
\end{equation}%
The segment $\mathcal{T}_{\left[ P_{0}P_{1}\right] }$ is cigar-shaped
two-dimensional surface. Let the length of the segment be $l\gg d$ and $\tau 
$, $0<\tau <l$ be a parameter along the segment. The radius $\rho $ of the
hollow segment tube as a function of $\tau $ has the form%
\begin{equation}
\rho ^{2}=\rho ^{2}\left( \tau \right) =\frac{1}{4}\frac{d}{\left(
l-d\right) ^{2}}\left( 2\tau -d\right) \left( 2l-3d\right) \left( 2l-2\tau
-d\right) ,\qquad 2d<\tau <l-2d  \label{a5.7}
\end{equation}%
If $\tau \gg 2d$, we obtain approximately 
\begin{equation}
\rho ^{2}=\rho ^{2}\left( \tau \right) =\frac{2d}{l}\tau \left( l-\tau
\right) ,\qquad 2d\ll \tau <l-2d  \label{a5.8}
\end{equation}%
The maximal radius of the segment tube $\rho _{\max }=\sqrt{ld/2}$ at $\tau
=l/2$.

It is clear that the ball $\mathcal{B}$ cannot be constructed only of such
tube segments. These hollow segments could not fill the ball completely. The
problem of construction of the ball, consisting of basic elements (points
and segments), appears to be a very complicated problem in the modified
geometry $\mathcal{G}_{\mathrm{d}}$. This problem of the ball construction
associates with the impossibility of deducing the geometry $\mathcal{G}_{%
\mathrm{d}}$ from a system of axioms. It seems that the geometry, deduced
from the system of axioms, cannot be multivariant, because any proposition,
obtained from axioms by means of the formal logic is to be definite. It
cannot contain different versions.

On the other hand, the multivariance is an essential property of the real
microcosm. It is a reason of quantum effects and atomism. The space-time
geometry is a basis of dynamics in microcosm.

\label{07b}In the multivariant geometry one constructs geometrical objects
by means of the deformation principle. Geometrical object $\mathcal{O}_{%
\mathrm{E}}$ is constructed in some region $\mathcal{S}_{1}$ of the proper
Euclidean space. In means that all blocks of $\mathcal{O}_{\mathrm{E}}$ as
well as the geometrical object $\mathcal{O}_{\mathrm{E}}$ itself are
expressed in terms of the world function $\sigma _{\mathrm{E}}$ of the
proper Euclidean geometry. Let us imagine that we need to shift $\mathcal{O}%
_{\mathrm{E}}$ in other region $\mathcal{S}_{2}$ of the space. We may move
all blocks of $\mathcal{O}_{\mathrm{E}}$ from $\mathcal{S}_{1}$ to $\mathcal{%
S}_{2}$ and construct a geometrical object $\mathcal{O}_{\mathrm{E}}^{\prime
}$ in $\mathcal{S}_{2}$ from blocks, using the same prescription of
construction, which has been used at construction of $\mathcal{O}_{\mathrm{E}%
}$. This prescription, written in terms of the world function, has the same
form in any geometry, if one uses only points as the basic concept of the
geometry. Points considered as blocks are not deformed at motion from the
region $\mathcal{S}_{1}$ to the region $\mathcal{S}_{2}$, even if the
geometry in the region $\mathcal{S}_{2}$ distinguishes from the geometry in
the region $\mathcal{S}_{1}$ (different world functions in the regions $%
\mathcal{S}_{1}$ and $\mathcal{S}_{2}$).

This consideration explains application of the deformation principle, but it
is not its proof. The deformation principle is the principle which admits
one to use physical geometries for description of the nonaxiomatizable
space-time geometry\label{07e}

\section{Concluding remarks}

Thus, there are three different representation of the proper Euclidean
geometry. In the E-representation there are three basic elements: point,
segment, angle. The segment and the angle are some auxiliary structures. The
segment is determined by two points. The angle is determined by three
points, or by two connected segments. In V-representation there are two
basic elements: point and vector (directed segment). The angle is replaced
by two segments (vectors). Its value is determined by the scalar product of
two vectors. Reduction of the number of basic elements is accompanied by
appearance of new structure: linear vector space with the scalar product,
given on it. Information, connected with the angle, is concentrated now in
the scalar product and in the linear vector space. Such a concept as
distance is a property of a vector (or a property of two points).

In the $\sigma $-representation there is only one basic element: point.
Interrelation of two points (segment) is described by the world function
(distance). Mutual directivity of segments, (or an angle) is considered as
an interrelation of three points. It is described by means of the scalar
product, expressed via distance (world function). In $\sigma $%
-representation the world function (distance) turns into a structure in the
sense, that the world function satisfies a series of constraints (conditions
I $\div $ IV). The concept of distance exists in all representations. But in
the E-representation and in the V-representation the distance is not
considered as a structure, because the conditions I $\div $ IV are not
considered as a constraints, imposed on the distance (world function). Of
course, these conditions are fulfilled in all representations, but they are
considered as direct properties of the proper Euclidean space, but not as
constraints, imposed on distance of the proper Euclidean space.

The $\sigma $-representation is interesting in the sense, that it contains
only one basic element (point) and only one structure (world function). All
other concepts of Euclidean geometry appear to be expressed via world
function. Supposing that these expressions have the same form in other
geometries, one can easily construct them, replacing world function. This
replacement looks as a deformation of the proper Euclidean space.

Usually a change of a representation is a formal operation, which is not
accompanied by a change of basic concepts. For instance, representations in
different coordinate systems differ only in the form of corresponding
algebraic expressions. A change of representation of the proper Euclidean
geometry is accompanied by a change of basic concepts, when the primary
concepts turns to the secondary ones and vice versa. \label{02b} The
procedure of a change of representation may be qualified as the logical
reloading. The logical reloading is rather rare logical procedure. Such a
change of primary concepts is unusual and difficult for a perception.

In particular, such difficulties of perception appear, because the linear
vector space is considered as an attribute of the geometry (but not as an
attribute of the geometry description). The linear vector space is an
attribute of the V-representation of the proper Euclidean geometry. There
are geometries (physical geometries), where the linear vector space cannot
be introduced at all. Physical geometries are described completely by the
world function, and the world function is the only characteristic of the
physical geometry. All other attributes of the physical geometry are
derivative. They can be introduced only via the world function. The physical
geometry cannot be axiomatized, in general. The proper Euclidean geometry $%
\mathcal{G}_{\mathrm{E}}$ is an unique known example of the physical
geometry, which can be axiomatized. Axiomatization of the proper Euclidean
geometry $\mathcal{G}_{\mathrm{E}}$ is used for construction of $\mathcal{G}%
_{\mathrm{E}}$. When the proper Euclidean geometry $\mathcal{G}_{\mathrm{E}}$
is constructed (deduced from the Euclidean axiomatics), one uses the fact
that the $\mathcal{G}_{\mathrm{E}}$ is a physical geometry. One expresses
all geometrical objects of the proper Euclidean geometry $\mathcal{G}_{%
\mathrm{E}}$ in terms of the Euclidean world function $\sigma _{\mathrm{E}}$%
. Replacing $\sigma _{\mathrm{E}}$ in all definitions of the $\mathcal{G}_{%
\mathrm{E}}$ by another world function $\sigma $, one obtains all
definitions of another physical geometry $\mathcal{G}$. It means that one
obtains another physical geometry $\mathcal{G}$, which cannot be axiomatized
(and deduced from some axiomatics)

Impossibility of the geometry $\mathcal{G}$ axiomatization is conditioned by
the fact, that the equivalence relation is intransitive, in general, in the
geometry $\mathcal{G}$. However, in any mathematical model, as well as in
any geometry, which can be axiomatized, the equivalence relation is to be
transitive. Almost all mathematicians believe, that any geometry can be
axiomatized. Collision of this belief with the physical geometry application
generates misunderstandings and conflicts \cite{R2008}.

\end{document}